\documentclass[english,journal]{IEEEtran}
\usepackage[T1]{fontenc}
\usepackage[latin9]{inputenc}
\usepackage{color}
\usepackage{babel}
\usepackage{amsmath}
\usepackage{amsthm}
\usepackage{amssymb}
\usepackage{graphicx}
\usepackage[unicode=true,pdfusetitle,
 bookmarks=true,bookmarksnumbered=false,bookmarksopen=false,
 breaklinks=false,pdfborder={0 0 0},pdfborderstyle={},backref=false,colorlinks=true]
 {hyperref}

\makeatletter
\theoremstyle{plain}

\theoremstyle{remark}

\theoremstyle{definition}

\theoremstyle{assumption}

\usepackage[ruled]{algorithm2e}
\newcommand{\nosemic}{\renewcommand{\@endalgocfline}{\relax}}
\newcommand{\dosemic}{\renewcommand{\@endalgocfline}{\algocf@endline}}
\usepackage{cite}
\usepackage[paper=letterpaper,top=0.75in,bottom=0.75in,right=0.75in,left=0.75in]{geometry}
\usepackage{subfigure}
\usepackage{enumerate}
\newtheorem{assumption}{Assumption}

\makeatother
\newtheorem{lemma}{Lemma}
\newtheorem{remark}{Remark}

\newtheorem{theorem}{Theorem}

\providecommand{\definitionname}{Definition}
\providecommand{\remarkname}{Remark}
\providecommand{\theoremname}{Theorem}

\begin{document}
\newgeometry{top=1in,bottom=0.75in,right=0.75in,left=0.75in}
\IEEEoverridecommandlockouts 
\title{Parameter Update Laws for Adaptive Control using Barrier Constraints}
\author{Ashwin P. Dani  \thanks{A. P. Dani is with the Department of Electrical and Computer Engineering at University of Connecticut, Storrs, CT 06269.}
}
\maketitle

\begin{abstract}
In this paper, constrained parameter update laws for adaptive control are developed using barrier constraints. An interpretation of the parameter update law from a constrained optimization problem, in which a regularized Barrier saddle function is formulated to incorporate parameter constraints using inverse and logarithmic barrier functions from interior-point methods. The resulting constrained update law is integrated with an adaptive trajectory tracking controller, enabling online learning of the unknown system model parameters. Forward invariance of the parameter estimate is established and Lyapunov stability of the closed-loop system with the constrained parameter update law is derived. The effectiveness of the proposed constrained adaptive control law is demonstrated through simulations, which validate its ability to maintain parameter estimates within prescribed bounds while ensuring convergence to the true parameter values and achieving steady state tracking performance.
\end{abstract}

\section{Introduction}

Adaptive control provides a framework to achieve control objectives in the presence of parametric uncertainties in the system model \cite{annaswamy2021historical, ioannou1996robust, narendra2012stable, krstic1995nonlinear, krause1990parameter, Slotine1991, volyanskyy2009new}. Various parameter update laws have been developed based on Lyapunov analysis and optimization techniques, with the gradient-based parameter update law being one of the most basic approaches for achieving regulation and tracking objectives. Several modifications of the gradient parameter update law exist, including the $\sigma$-modification \cite{ioannou1996robust}, the e-modification \cite{narendra2012stable}, and $\mathcal{L}_1$-norm-based parameter update laws \cite{cao2008design} for constant parameters, as well as update laws in \cite{chen2020adaptive, gaudio2021parameter} for time-varying parameters. These parameter update laws typically require the persistence of excitation (PE) condition to be satisfied for parameter convergence \cite{chowdhary2013concurrent}, which is difficult to verify in real time. Parameter update laws using regressor filtering methods are developed in \cite{kreisselmeier1982rate} based on interval excitation (IE), in \cite{roy2017combined} using initial excitation, and more recent results in \cite{ortega2020new} demonstrate parameter convergence under conditions weaker than PE. Building on least-squares concepts, concurrent learning (CL)-based parameter update laws are proposed in \cite{chowdhary2013concurrent}, and integral concurrent learning (ICL) is introduced in \cite{parikh2019integral} to improve parameter convergence under a finite excitation (FE) condition, which is a weaker condition compared to PE and IE. A parameter estimator achieving finite time parameter identification is developed in \cite{adetola2008finite}.
More recently, momentum-based tools from optimization theory have been used to design parameter update laws in the context of reinforcement learning controllers \cite{somers2024online}, with the goal of improving transient performance. One of the primary objectives of these recent parameter update law designs is to achieve convergence of the parameter estimates to their true values, thereby improving robustness and transient performance while relaxing the PE condition \cite{parikh2019integral}.

In adaptive control, it is well known that parameters estimated using gradient-based update laws can drift \cite{anderson1985adaptive, ioannou1996robust}. To address this issue, several modifications have been proposed, including the $\sigma$-modification \cite{ioannou1996robust}, the e-modification \cite{narendra2012stable}, and smooth projection-based methods \cite{ioannou1996robust}. A simpler projection-like method for adaptive control is derived in \cite{bakker1996stability}. These update laws typically guarantee only that the parameter estimates remain bounded within a compact neighborhood of a preselected region \cite{chowdhary2013concurrent}. The motivation for imposing constraints on the parameter update laws is to accelerate convergence and reduce large transients \cite{ioannou1996robust}. A method developed in \cite{akella2005novel} guarantees strict boundedness of the parameter estimates within prescribed lower and upper bounds. In certain robotics and space applications, additional constraints on the parameters, such as symmetry with eigenvalue bounds or positive definiteness of parameter matrices, are important to ensure physically consistent parameter estimates and to improve controller robustness. Methods in \cite{boffi2021implicit, lee2018natural} develop parameter update laws based on Bregman divergence metrics to enforce positive definiteness of parameter matrices. The method in \cite{moghe2022projection} incorporates symmetric matrix constraints with eigenvalue bounds into the parameter estimation law. These approaches allow the incorporation of specific structural constraints present in the system dynamics. In some applications, it is critical to strictly enforce bound constraints on the parameters. 

Studies in \cite{boffi2021implicit,annaswamy2021historical} highlight that gradient update laws can also be derived using speed-gradient methods, as shown in \cite{fradkov1979speed,fradkov2013nonlinear}. In these approaches, an objective function is explicitly derived, from which a gradient-based update law is obtained by taking the negative gradient of the function. In our recent study \cite{DaniAACRL2025}, a constrained parameter update law based on the inverse barrier method for reinforcement learning is developed, which incorporates upper and lower bounds on each element of the parameter vector when uncertainty is present in the control effectiveness matrix. The generalization of the inverse barrier-based parameter update law to the adaptive control setting is presented by considering additional constraints, such as norm constraints on the estimated parameters \cite{DaniCDC2025}.

The contribution of this paper is the design of a new constrained parameter update law for adaptive control using barrier function method. Constrained parameter update laws based on both inverse and logarithmic barrier functions are developed. To design the constrained parameter update law, a regularized Barrier saddle (RBS) function is formulated using an objective function that is convex in the parameter estimates, and strictly concave in dynamic Barrier weights. To enforce the parameter constraints, inverse and logarithmic barrier methods are employed \cite{Boyd2004,bertsekas2014constrained}. The applicability of the proposed constrained parameter update law to incorporate upper and lower bounds on the components of the parameter estimates, is demonstrated, thereby illustrating the generality of the approach. A Lyapunov-based stability analysis for the adaptive controller with the constrained parameter update law is developed. The Lyapunov analysis shows that the closed-loop system is uniformly ultimately bounded (UUB) for a trajectory tracking objective, while the parameter estimates remain strictly within the prescribed constraints, provided they are initialized within the bounds. Simulation studies demonstrate that the proposed constrained parameter update law ensures parameter constraints are satisfied. In contrast, gradient-based and CL-based parameter update laws may not always guarantee that the parameter constraints are satisfied.
\noindent\paragraph*{Notations} For a vector constraint $c_{j}(\theta): \mathbb{R}^p \rightarrow \mathbb{R}^p$, $j \in \{1,2\}$ with components $c_{ji}$ such that $c_j= [c_{j1},...,c_{jp}]^T$ with $\theta \in \mathbb{R}^p$, $\nabla_\theta c_j$ denotes component-wise gradient defined by $\nabla_{\hat{\theta}} c_j = [\frac{\partial c_{j1}}{\partial \hat{\theta}_1},...,\frac{\partial c_{jp}} {\partial \hat{\theta}_p}]^T$. The projection operator $[a]_b^+$ for $b\in \mathbb{R}_{\geq 0}$ is defined as 
\begin{equation}
    [a]_b^+ = \begin{cases}
      a, & \text{if           $b>0$}\\
      \mathrm{max}\{0,a\}, & \text{if           $b=0$}
    \end{cases}
\end{equation}
For a differentiable convex function $c_{ij}(\hat{\theta})$ in $\hat{\theta}$, the inequality is satisfied $c_{ji}(\theta) \geq c_{ji}(\hat{\theta}) + \nabla_{\hat{\theta}} c_{ij}(\hat{\theta})^T(\theta - \hat{\theta})$. The symbols $\mathbb{R}_{\geq 0} = [0,\infty)$, $\mathbb{R}_+ = (0,\infty)$.

\section{System Model and Control Objective}
\subsection{System Dynamics}
Consider the following system
\begin{equation}
\dot{x}=f(x) + u \label{eq:SystemModel}
\end{equation}
where $x(t) \in \mathbb{R}^n$ is the system state, $u\in\mathbb{R}^n$ is the control input, $f:\mathbb{R}^n\rightarrow \mathbb{R}^n$ is a locally Lipschitz continuous function which is linearly parametrizable, defined as
\begin{align}
f(x) = Y(x)\theta \label{eq:LinParam}
\end{align}
where $Y: \mathbb{R}^n \rightarrow \mathbb{R}^{n \times p}$ is the regressor matrix, which is continuous and locally bounded in $x$, and $\theta \in \mathbb{R}^p$ is an unknown parameter.

\subsection{Controller Objective}
The control objective is to track a desired trajectory $x_d(t)$ with properties that $x_d \in \mathcal{C}^1$, and $x_d(t),\dot{x}_d(t) \in \mathcal{L}_\infty$ and compute parameter vector estimate $\hat{\theta}(t)\in \mathbb{R}^p$ while maintaining prescribed convex constraints on the parameter estimates, such as upper and lower bounds defined as $\underline{\theta}_i < \hat{\theta}_i < \bar{\theta}_i$ or convex norm constraints of the form $ \Vert \hat{\theta} \Vert < \bar{\theta}$.
For the control design, let's define tracking and parameter estimation errors as follows
\begin{equation}
    e(t) = x(t) - x_d(t), \quad
    \tilde{\theta}(t) = \theta - \hat{\theta}(t)
    \label{eq:regulation_error}
\end{equation}
where $\tilde{{\theta}}(t)\in\mathbb{R}^p$ is the parameter estimation error.

\subsection{Adaptive Controller Design} \label{sec:ControlDesign}
To achieve the control objective, an adaptive trajectory tracking controller is designed as
\begin{equation}
u = \dot{x}_d - Y\hat{\theta} - ke \label{eq:Control}
\end{equation}
where $\dot{x}_d(t)$ is the derivative of the desired trajectory, and $k>0$ is the control gain. Taking the time derivative of (\ref{eq:regulation_error}) and substituting the controller (\ref{eq:Control}), the following closed-loop error system can be obtained 
\begin{equation}
    \dot{e} = Y\tilde{\theta} -ke  \label{eq:ErrorDyn}
\end{equation}
The parameter estimation error dynamics is written as
\begin{equation}
\dot{\tilde{\theta}} = -\dot{\hat{\theta}}
\end{equation}

\subsection{Design of constrained adaptive parameter update law}
In this section, the inverse barrier method and the logarithmic barrier method are used to incorporate constraints on the parameter estimates. The constrained parameter update law is first presented using upper and lower bound constraints on the components of the parameter estimates. 
\begin{assumption}
    The components of the true parameters are bounded by $\underline{\theta}_i + \delta_{1i} < \theta_i < \bar{\theta}_i - \delta_{2i}, \quad \forall i=\{1,...,p\}$, where $\delta_{1i}$ and $\delta_{2i}$ are positive margin parameters.    \label{ass:MarginParam}
\end{assumption}
Let $\Theta_c \subset \mathbb{R}^p$ be an open and bounded constrained set. Consider the following parameter constraint set
\begin{equation}
    \Theta^c = \{\hat{\theta}: \underline{\theta}_i < \hat{\theta}_i(t) < \bar{\theta}_i, \;\; \forall i=\{1,...,p\}\} \label{eq:ParamInterior}
\end{equation}
which can be written as
\begin{equation}
    g_{1i} = \underline{\theta}_i - \hat{\theta}_i(t) < 0, \quad g_{2i} = \hat{\theta}_i(t) - \bar{\theta}_i < 0.
\end{equation}
Using inverse Barrier method, the constraints with scaling and a shift are formulated as
\begin{equation}
    c_{1i}(\hat{\theta}) = \frac{-\delta_{1i}}{\underline{\theta}_i- \hat{\theta}_i(t)}-1, \quad c_{2i}(\hat{\theta}) = \frac{-\delta_{2i}}{\hat{\theta}_i(t) - \bar{\theta}_i} - 1 \label{eq:InvBarConst}
\end{equation}
Similarly, using log-Barrier method, the constraints with scaling are formulated as
\begin{equation}
    c_{1i}(\hat{\theta}) =-\mathrm{log}\left(\frac{-\underline{\theta}_i+ \hat{\theta}_i(t)}{\delta_{1i}}\right), c_{2i}(\hat{\theta}) = -\mathrm{log}\left(\frac{-\hat{\theta}_i(t) + \bar{\theta}_i}{\delta_{2i}}\right). \label{eq:logBarConst}
\end{equation}
For CL-based parameter update law, consider a history stack $\mathcal{H} = \{x_k,\hat{u}_k,\dot{\hat{x}}_k\}_{k=1}^{k=m}$ which starts collecting data as the system evolves. Let $S(t) = \sum_{k=1}^m Y_k^TY_k$ be a matrix computed using history stack data, and there exists finite $t_H$ such that $t_H = \mathrm{inf}\{t \geq t_0: \lambda_\mathrm{min}(S(t)) \geq \bar{\sigma}_1 \}$. Let $\xi(t)$ be a switching signal such that $\xi(t) = 0$ for $t<t_H$ and $\xi(t) = 1$ for $t\geq t_H$. 

The constrained parameter update law can be written in the following form
\begin{align}
\dot{\hat{\theta}} =& P Y(x)^T e + Pk_{cls} \sum_{k=1}^m Y_k^T(\dot{\hat{x}}_k - u_k - Y_k\hat{\theta}(t)) \nonumber \\
&- \sum_{j=1}^2 P \mathrm{diag}(\lambda_j)\nabla_{\hat{\theta}}c_j \nonumber \\
\dot{\lambda}_j =& [-\alpha \lambda_j + \Gamma_j^{-1} c_j]^+_{\lambda_j}, \quad \forall j = \{1,2\}  \label{eq:ParameterUpdateLaw}
\end{align}
where $P \in \mathbb{R}^{p \times p}$, $\Gamma_j \in \mathbb{R}^{p \times p}$ are positive definite (PD) diagonal learning rate matrices, $\lambda = [\lambda_{1}^T, \lambda_{2}^T]^T \in \mathbb{R}^{2p}$ are non-negative dynamic Barrier weights such that $\lambda_j(t_0)>0_{p \times 1}$, $\alpha >0$ and $k_{cls}=\xi(t)k_{cl}$ where $k_{cl}>0$ are constant gains. The gradient $\nabla_{\hat{\theta}} c_j$ is computed component-wise, i.e., $\nabla_{\hat{\theta}} c_j = [\frac{\partial c_{j1}}{\partial \hat{\theta}_1},...,\frac{\partial c_{jp}} {\partial \hat{\theta}_p}]^T$. 
\begin{remark}
    The margin parameters $\delta_{1i}$ and $\delta_{2i}$ define interior thresholds relative to the lower and upper parameter estimate bounds, respectively. For both the shifted inverse-barrier and scaled log-barrier formulations in (\ref{eq:InvBarConst})-(\ref{eq:logBarConst}), $c_{1i}(\hat\theta)=0$ when $\hat\theta_i-\underline{\theta}_i=\delta_{1i}$, and $c_{2i}(\hat\theta)=0$ when $\bar{\theta}_i-\hat\theta_i=\delta_{2i}$. Consequently, Assumption 1 ensures $c_{1i}(\theta)\leq0$ and $c_{2i}(\theta)\leq0$. The margin parameters do not alter the feasible set $\Theta_c$, rather, they determine an interior region in which the barrier functions are nonpositive and regulate how early the dynamic barrier weights respond as the parameter estimates approach the constraint boundaries.
\end{remark}
\begin{remark}
    The constrained parameter update law is designed based on CL parameter update law so that the PE condition for convergence of the parameter is not required for parameter convergence. 
\end{remark}
\begin{remark}
    The parameter estimates are initialized such that the constraints are satisfied at $t=t_0$.
\end{remark}

\begin{assumption} 
For the history stack $\mathcal{H} = \{x_k,\hat{u}_k,\dot{\hat{x}}_k\}_{k=1}^{k=m}$ the following condition is satisfied
\begin{equation}
    \lambda_{min}(S(t)) = \lambda_{min}(\sum_{k=1}^m Y_k^T Y_k) \geq \bar{\sigma}_1 > 0, \;\;\forall t \geq t_H 
\end{equation}
where $\bar{\sigma}_1 \in \mathbb{R}^+$. The numerically computed derivatives of $x(t)$, $\dot{\hat{x}}_k$ computed at $k$th data point satisfies $\dot{\hat{x}}_k - \dot{x}_k = \epsilon_{xk}$ where $\epsilon_{xk} \in \mathbb{R}^n$ such that $\Vert \epsilon_{xk} \Vert \leq \bar{\epsilon}_x$. \label{ass:FiniteExcitation}
\end{assumption}
\begin{remark}
Assumption \ref{ass:FiniteExcitation} is a finite excitation condition that can be verified in real-time \cite{chowdhary2013concurrent}.
\end{remark}
\begin{assumption} \label{ass:HistoryStackRegularity}
    The symmetric positive semidefinite matrix $S(t)$ is piecewise continuous and uniformly bounded in $t$. Moreover, the history stack update times are locally finite. 
\end{assumption}
\subsubsection{Gradients of Inverse and Log-barrier Constraints}
For the parameter update law in (\ref{eq:ParameterUpdateLaw}), the gradients of inverse and log-barrier are shown below. Compared to the inverse Barrier, log-Barrier constraints generally show better numerical stability near the constraints\cite{Boyd2004, bertsekas2014constrained}. The gradients of inverse-Barrier constraints $c_1(\hat{\theta})$ and $c_2(\hat{\theta})$ are given by
\begin{align}
    \nabla_{\hat{\theta}}c_1 &= \left[\frac{-\delta_{11}}{(-\underline{\theta}_1 + \hat{\theta}_1(t))^2}, ..., \frac{-\delta_{1p}}{(-\underline{\theta}_p + \hat{\theta}_p(t))^2}\right]^T  \nonumber \\
    \nabla_{\hat{\theta}}c_2 &= \left[\frac{\delta_{21}}{(- \hat{\theta}_1(t) + \bar{\theta}_1)^2}, ..., \frac{\delta_{2p}}{(- \hat{\theta}_p(t) + \bar{\theta}_p)^2}\right]^T 
\end{align}
The gradients of log-Barrier constraints $c_1(\hat{\theta})$ and $c_2(\hat{\theta})$ are given by
\begin{align}
    \nabla_{\hat{\theta}}c_1 &= \left[\frac{-1}{-\underline{\theta}_{1} + \hat{\theta}_1(t)}, ..., \frac{-1}{-\underline{\theta}_{p} + \hat{\theta}_p(t)}\right]^T  \nonumber \\
    \nabla_{\hat{\theta}}c_2 &= \left[\frac{1}{- \hat{\theta}_1(t) + \bar{\theta}_{1}}, ..., \frac{1}{- \hat{\theta}_p(t) + \bar{\theta}_{p}}\right]^T 
\end{align}
\subsubsection{Derivation of constrained adaptive update law}
The Barrier constrained parameter update law can be intuitively derived by solving the following constrained minimization problem
\begin{align}
    &\hat{\theta}^* = \mathrm{min}_{\hat{\theta}} \; e(t)^TY(t)\tilde{\theta} + \frac{k_{cls}}{2} \left[\sum_{k=1}^m \Vert \hat{\dot{x}}_k - u_k -Y_k\hat{\theta}\Vert^2\right]  \nonumber \\
    \mathrm{s.\;t.}&  \; g_{1i}(\hat{\theta}) <0, \quad g_{2i}(\hat{\theta}) <0, \;\; \forall \; i = \{1,...,p\}. \label{eq:unconProb}
\end{align}
Motivated by the constrained Barrier (inverse or log) problem in (\ref{eq:unconProb}), let's define the following time-varying convex-concave regularized Barrier saddle (RBS) function $L(\hat{\theta},\lambda,t): \mathbb{R}^p \times \mathbb{R}^{2p} \times \mathbb{R}^+ \rightarrow\mathbb{R}$, defined by
\begin{align}
    L(\hat{\theta},\lambda,t) &=  e(t)^TY(t)\tilde{\theta} + \frac{k_{cls}}{2} \left[\sum_{k=1}^m \Vert \hat{\dot{x}}_k - u_k -Y_k\hat{\theta}\Vert^2\right] \nonumber \\
    &-\frac{\alpha}{2}\lambda^T\Gamma\lambda + \sum_{j=1}^2  \lambda_j^T c_j \label{eq:Lagrangian}
\end{align}
Note that RBS function $L$ is convex in $\hat{\theta}$ and strongly concave in $\lambda$. After Assumption \ref{ass:FiniteExcitation} is satisfied, the CL quadratic term makes $\mathcal{L}$ strongly convex in $\hat{\theta}(t)$. The parameter update law in (\ref{eq:ParameterUpdateLaw}) can be computed using $\dot{\hat{\theta}} = -P\nabla_{\hat{\theta}}L$.
The $\lambda$ update rule of (\ref{eq:ParameterUpdateLaw}) can be developed using the strong concavity of $L$ in $\lambda$, and $\dot{\lambda}_j = [\Gamma^{-1}\nabla_{\lambda}L]^+_{\lambda_j} $.
Due to projection law and switching of gain $k_{cls}$ based on history stack rank condition, the solution of (\ref{eq:ParameterUpdateLaw}) is interpreted in the sense of Caratheodory, which is differentiable at almost all $t \in \mathbb{R}^+$, starting from the initial condition $\hat{\theta}(t) \in \Theta^c$ and  $\lambda_j(t_0)$ > 0.

The error dynamics of the parameter update law can now be written as
\begin{align}
    \dot{\tilde{\theta}} &= -PY^{T}e - Pk_{cls} \sum_{k=1}^m Y_k^TY_k\tilde{\theta} \nonumber \\
    &+ \sum_{j=1}^2 P\mathrm{diag}(\lambda_j)\nabla_{\hat{\theta}} c_j - Pk_{cls}d_x  \label{eq:thetaTildeDot}
\end{align}
where $d_x =\sum_{k=1}^mY_k^T\epsilon_{xk}$.

\section{Forward Invariance and Stability Analysis}
In this section, the forward invariance of the parameter update law is analyzed, which will establish that $\hat{\theta} \in \Theta^c$. Consider a state $\zeta = [e^T,\hat{\theta}^T,\lambda^T]^T$ and a domain $\mathcal{D} =\mathbb{R}^n \times \Theta^c \times \mathbb{R}^{2p}_{\geq 0}$. The system (\ref{eq:ErrorDyn})-(\ref{eq:ParameterUpdateLaw}) can be represented as a time-varying projected dynamical system, whose unprojected right hand side is continuous in $\zeta$ and $S(t)$ is piecewise continuous using Assumption \ref{ass:HistoryStackRegularity}. Under these regularity assumptions, the existence of a Caratheodory solution $\zeta(t)$ is ensured locally on a maximal time interval $[t_0,\tau_\mathrm{max})$ for $\zeta(t_0) \in \mathcal{D}$ \cite{heemels2023existence}. Let $\zeta : [t_0,\tau_{\mathrm{max}}) \rightarrow \mathcal{D}$ be a maximal Caratheodory solution of (\ref{eq:ErrorDyn}), (\ref{eq:ParameterUpdateLaw}), initialized at $\zeta(t_0) \in \mathcal{D}$. We prove that $\tau_\mathrm{max} = \infty$, which proves the forward invariance of $\hat{\theta}(t)$ within $\Theta^c$.

Let's define distance from the boundaries $d_{1i} = \hat{\theta}_i - \underline{\theta}_i$ and $d_{2i} = \bar{\theta}_i - \hat{\theta}_i$ and $\beta_{ji}:\mathbb{R}_+ \rightarrow \mathbb{R}_+$ defined as
\begin{equation}
    \beta_{ji}(d_{ji}) = 
    \begin{cases}
    \delta_{ji}/d_{ji}^2, \quad \text{shifted inverse barrier}, \\
    1/d_{ji}, \quad \text{logarithmic barrier}.
\end{cases}  \label{eq:betaDef}
\end{equation}
where $\beta_{ji} \rightarrow \infty$ as $d_{ji} \rightarrow 0^+$. Then the gradients $\frac{\partial c_{1i}}{\partial \hat{\theta}_i} = - \beta_{1i}(d_{1i})$ and $\frac{\partial c_{2i}}{\partial \hat{\theta}_i} =  \beta_{2i}(d_{2i})$. The dynamics of $d_{ji}$ are given by
\begin{align}
    \dot{d}_{1i} = - \dot{d}_{2i} &= \dot{\hat{\theta}}_i = q_i
    + p_{ii}\lambda_{1i}\beta_{1i}(d_{1i}) - p_{ii}\lambda_{2i}\beta_{2i}(d_{2i}) \label{eq:GapDot}
\end{align}
where $q_i = p_{ii}[Y^Te + k_{cls}S(t)\tilde{\theta}(t)+k_{cls}d_x]_i$, and $S(t) = \sum_{k=1}^m Y_k^TY_k$.
\begin{lemma} \label{lem:BoundednessLemma}
    For a maximal Caratheodory solution of the system (\ref{eq:ErrorDyn}), (\ref{eq:ParameterUpdateLaw}), provided that Assumptions \ref{ass:MarginParam}-\ref{ass:HistoryStackRegularity} are satisfied, the closed-loop signals are bounded on $t \in [t_0,\tau_{\mathrm{max}})$ as follows
    \begin{equation}
        \Vert e(t)\Vert \leq e_m, \;  \Vert \tilde{\theta}(t)\Vert \leq \tilde{\theta}_m,\; \Vert \lambda(t)\Vert \leq \lambda_m \label{eq:Lemma1Bounds}
    \end{equation}
    where $e_m$, $\tilde{\theta}_m$, and $\lambda_m$ are positive  constants.
\end{lemma}
\begin{proof}
Since $\lambda_{ji}(t_0) \geq 0$, from the definition of projection, at $\lambda_{ji} = 0$, $\dot{\lambda}_{ji} = \mathrm{max}\{0,-\alpha \lambda_{ji} + (\Gamma_j^{-1})_{ii} c_{ji}\} \geq 0$. Thus, $\lambda_{ji}(t) \geq 0, \; \forall t< \tau_\mathrm{max}$.

    Let $z(t) = [e(t)^T\; \tilde{\theta}(t)^T \; \lambda(t)^T]^T \in \mathbb{R}^{n+3p}$ be an auxiliary vector of closed-loop signals. Consider a following Lyapunov function candidate, $V: \mathbb{R}^n \times \mathbb{R}^{p} \times \mathbb{R}^{2p}_{\geq 0} \rightarrow \mathbb{R}_{\geq 0}$
\begin{equation}
    V(z) = \frac{1}{2}e^Te + \frac{1}{2}\tilde{\theta}^TP^{-1}\tilde{\theta} + \frac{1}{2}\lambda^T\Gamma\lambda \label{eq:LyapunovFun}
\end{equation}
where $\Gamma = \mathrm{blkdiag}(\Gamma_1,\Gamma_2)$, and satisfies 
\begin{equation}
\Lambda_{\mathrm{min}}\Vert z \Vert ^2 \leq V(z) \leq \Lambda_{\mathrm{max}}\Vert z \Vert ^2    \label{eq:VBound}
\end{equation}
where $\Lambda_{\mathrm{min}} = \mathrm{min}\{\frac{1}{2},\frac{1}{2}\lambda_{\mathrm{min}}(P^{-1}),\frac{1}{2}\lambda_{\mathrm{min}}(\Gamma)\}$, $\Lambda_{\mathrm{max}} = \mathrm{max}\{\frac{1}{2},\frac{1}{2}\lambda_{\mathrm{max}}(P^{-1}),\frac{1}{2}\lambda_{\mathrm{max}}(\Gamma)\}$.
Taking the time derivative of $V(z)$ yields
\begin{align}
\dot{V} &= e^T\dot{e} + \tilde{\theta}^TP^{-1}\dot{\tilde{\theta}} + \lambda^T\Gamma \dot{\lambda} \nonumber \\
\dot{V} &= e^TY\tilde{\theta} - ke^Te + \tilde{\theta}^TP^{-1}(-PY^Te- Pk_{cls} \sum_{k=1}^m Y_k^TY_k\tilde{\theta} \nonumber \\
&+ \sum_{j=1}^2 P\mathrm{diag}(\lambda_j)\nabla_{\hat{\theta}} c_j -Pk_{cls}d_x) \nonumber \\
&+ \sum_{j=1}^2\lambda_j^T\Gamma_j[-\alpha \lambda_j + \Gamma_j^{-1}c_j]^+_{\lambda_j} \nonumber \\
\dot{V} &= -ke^Te - k_{cls} \tilde{\theta}^T\sum_{k=1}^m Y_k^TY_k\tilde{\theta} + \tilde{\theta}^T(\sum_{j=1}^2 \mathrm{diag}(\lambda_j)\nabla_{\hat{\theta}} c_j) \nonumber \\
&+ \sum_{j=1}^2 \lambda_j^T\Gamma_j(-\alpha \lambda_j + \Gamma_j^{-1}c_j) - k_{cls}\tilde{\theta}^T d_x
\end{align}
where $\lambda_{ji}\Gamma_{jj}[-\alpha \lambda_{ji} + \Gamma_{jj}^{-1}c_{ji}]^+_{\lambda_{ji}} = \lambda_{ji}\Gamma_{jj}(-\alpha \lambda_{ji} + \Gamma_{jj}^{-1}c_{ji})$ since for $\lambda_{ji} > 0$ the projection is inactive and for $\lambda_{ji} = 0$, both sides are $0$.
The following bound on $\dot{V}(z)$ can now be obtained
\begin{align}
\dot{V} &\leq - ke^Te - k_{cls} \tilde{\theta}^T\sum_{k=1}^m Y_k^TY_k\tilde{\theta}  - \alpha \underline{\gamma} \sum_{j=1}^2 \lambda_j^T \lambda_j \nonumber \\
&+ \tilde{\theta}^T(\sum_{j=1}^2 \mathrm{diag}(\lambda_j)\nabla_{\hat{\theta}} c_j) + \sum_{j=1}^2 \lambda_j^Tc_j - k_{cls}\tilde{\theta}^T d_x   \label{eq:VDeriv}
\end{align}
where $\underline{\gamma} = \lambda_{min}(\Gamma)$. 
Using the convexity of $c_{ji}(\hat{\theta})$ in $\hat{\theta}$
\begin{equation}
    c_{ji}(\theta) \geq c_{ji}(\hat{\theta}) + \nabla_{\hat{\theta}} c_{ji}(\hat{\theta})^T(\theta - \hat{\theta}) 
\end{equation}
where $\nabla_{\hat{\theta}} c_{ji}$ is a standard gradient. Multiplying by $\lambda_{ji} \geq 0$, for $j=\{1,2 \}$ and summing over $i=\{1,..,p \}$ yields
\begin{equation}
    \tilde{\theta}^T \mathrm{diag}(\lambda_j)\nabla_{\hat{\theta}} c_j + \lambda_j^T c_j(\hat{\theta}) \leq \lambda_j^T c_j(\theta) \label{eq:convexInequal}
\end{equation}
Using (\ref{eq:convexInequal}) results in the following bound on $\dot{V}(z)$
\begin{align}
    \dot{V} &\leq - ke^Te - k_{cls} \tilde{\theta}^T\sum_{k=1}^m Y_k^T Y_k\tilde{\theta}  - \alpha \underline{\gamma} \sum_{j=1}^2 \lambda_j^T \lambda_j \nonumber \\
    &- k_{cls}\tilde{\theta}^T d_x+ \sum_{j=1}^2 \lambda_j^Tc_j(\theta) \label{eq:VDerviInter}
\end{align}
The Assumption \ref{ass:MarginParam} follows
\begin{equation}
    c_{1i}(\theta) \leq 0, \;\;\; c_{2i}(\theta) \leq 0 \label{eq:constraintInequLem1}
\end{equation}
for both inverse and log Barrier constraints in (\ref{eq:InvBarConst}) and (\ref{eq:logBarConst}), and given that the projection ensures $\lambda_{ji} \geq 0$, leads to $\sum_{j=1}^2 \lambda_j^Tc_j(\theta) \leq 0$. Since $k_{cls}=0$ when $t<t_H$, we have
thus, a.e. on $[t_0,\mathrm{min}(t_H,\tau_\mathrm{max}))$
\begin{align}
\dot{V} &\leq - k\Vert e \Vert^2 -\alpha\underline{\gamma}\Vert \lambda \Vert^2  \leq 0 \label{eq:VdotBound1Lem1}
\end{align}
Hence, if $t_H \leq \tau_\mathrm{max}$, $V(t_H)\leq V(t_0)$. For $t\geq t_H$, we have
\begin{align}
\dot{V} &\leq - k\Vert e \Vert^2 - k_{cl}\bar{\sigma}_1\Vert \tilde{\theta} \Vert^2  -\alpha\underline{\gamma}\Vert \lambda \Vert^2 + k_{cl}\Vert \tilde{\theta} \Vert \bar{d}_x  \nonumber \\
&\leq - k\Vert e \Vert^2 - \frac{k_{cl}\bar{\sigma}_1}{2}\Vert \tilde{\theta} \Vert^2  -\alpha\underline{\gamma}\Vert \lambda \Vert^2 + \frac{k_{cl}}{2\bar{\sigma}_1} \bar{d}_x ^2\label{eq:VdotBound2Lem1}
\end{align}
where $\Vert d_x \Vert \leq\bar{d}_x$. 
Using (\ref{eq:VBound}) and (\ref{eq:VdotBound2Lem1}), a.e. for $t \in [t_H,\tau_\mathrm{max})$, we have $\dot{V} \leq -\eta\Vert z \Vert^2 + D_x$, where $\eta = \mathrm{min}\{k,\frac{k_{cl}\bar{\sigma}_1}{2},\alpha \underline{\gamma} \}$ and $D_x =\frac{k_{cl}}{2\bar{\sigma}_1} \bar{d}_x ^2$. Using Comparison lemma, $V(t) \leq V(t_H)e^{-\beta(t-t_H)} + \frac{D_x}{\beta}(1-e^{-\beta(t-t_H)})$, where $\beta = \frac{\eta}{\Lambda_\mathrm{max}}$. From (\ref{eq:VdotBound1Lem1}), $V(t_H)\leq V(t_0)$, Therefore,
\begin{equation}
    V(t) \leq V_\mathrm{max} := \mathrm{max}\{V(t_0),\frac{D_x}{\beta} \}, \;\; t \in [t_0,\tau_{\mathrm{max}}).
\end{equation}
The bounds (\ref{eq:Lemma1Bounds}) follow from (\ref{eq:LyapunovFun}) and (\ref{eq:VBound}).
\end{proof}

\begin{lemma} (Forward Invariance)
    Provided $\zeta(t_0)\in\mathcal{D}$, Assumptions \ref{ass:MarginParam}-\ref{ass:HistoryStackRegularity} are satisfied, for every maximal Caratheodory solution of (\ref{eq:ErrorDyn})-(\ref{eq:ParameterUpdateLaw}), $\exists\epsilon_{1i},\epsilon_{2i}>0, \; i=\{1,..,p\}$ such that
    \begin{equation}
         \underline{\theta}_i + \epsilon_{1i}\leq \hat{\theta}(t) \leq \hat{\theta}_i - \epsilon_{2i},\; t \in [t_0,\tau_{\mathrm{max}})
    \end{equation}
    Furthermore, $\tau_\mathrm{max} = \infty$. Consequently, $\Theta^c$ is forward invariant.\label{lem:ForwardInvariance}
\end{lemma}
\begin{proof}

Consider a maximal Caratheodory solution $\zeta(t)$, $t\in[t_0,\tau_\mathrm{max})$. From Lemma  \ref{lem:BoundednessLemma}, (\ref{eq:Lemma1Bounds}) is satisfied for $t<\tau_\mathrm{max}$. Since $x_d(t) \in \mathcal{L}_\infty$, $x(t) \in \mathcal{L}_\infty$. Since $Y(x)$ is continuous, and locally bounded in $x$, we have $\Vert Y(x) \Vert \leq \bar{Y}$ and using Assumption \ref{ass:HistoryStackRegularity}, $\Vert S(t) \Vert \leq \bar{S}$, which yields following bound on $q_i$ of (\ref{eq:GapDot})
    \begin{equation}
        \vert q_i \vert \leq p_{ii}\bar{Y}e_m + p_{ii}k_{cl}\bar{S}\tilde{\theta}_m +p_{ii} k_{cl}\bar{d}_x \leq \bar{q} \label{eq:qBound}
    \end{equation}
    \paragraph{Behavior of $d_{1i}$ close to $\underline{\theta}_i$} In (\ref{eq:GapDot}), there are three terms on the right hand side. The first term $\vert q_i \vert \leq \bar{q}$. Consider the last term $p_{ii}\lambda_{2i}\beta_{2i}(d_{2i})$, where using the result of Lemma \ref{lem:BoundednessLemma}, $\lambda_{2i} \leq \lambda_m$.

    Now consider intervals near the lower boundary $0<d_{1i}<r_{1i}$, where $0<r_{1i}<\mathrm{min}\{d_{1i}(t_0), \delta_{1i},\frac{w_i}{2} \}$, where $w_i = d_{1i}+d_{2i}=\bar{\theta}_i - \underline{\theta}_i$. When $d_{li} \in (0,r_{1i})$, we have $\beta_{2i}(d_{2i}) \leq \beta(w_i - r_{1i}) = \bar{\beta}_2$, where $\bar{\beta}_2 \in \mathbb{R}^+$. Using the result of Lemma \ref{lem:BoundednessLemma}, $\lambda_{2i} \leq \lambda_m$. Then the last term $p_{ii}\lambda_{2i}\beta_{2i}(d_{2i})$ and first term can be bounded together by $M_{1i} >0$
    \begin{equation}
        \vert q_i - p_{ii}\lambda_{2i}\beta_{2i}(d_{2i}) \vert \leq M_{1i}, \; \mathrm{for } \;0<d_{1i}<r_{1i}
    \end{equation}
    Given $\lambda_{1i} \geq 0$ and $\beta_{1i} \geq 0$, this gives a lower bound on $\dot{d}_{1i}$ as
    \begin{equation}
        \dot{d}_{1i} \geq -M_{1i} + p_{ii} \lambda_{1i}\beta_{1i}(d_{1i}) \geq -M_{1i} \label{eq:dDotLowerBound}
    \end{equation}
    Now, it is checked if the positive second term of (\ref{eq:dDotLowerBound}) dominates the first negative term as $\hat{\theta}_i \rightarrow \underline{\theta}_i$, i.e., when $d_{1i} \rightarrow0^+$. For that we first need to check whether $\lambda_{1i}$ becomes non-zero positive as $d_{1i} \rightarrow 0^+$. 
    
    Consider an interval $[\frac{r_{1i}}{2},r_{1i}]$ and two time instances $t_1$ and $t_2$. Let $t_2>t_1$ be the first time $d_{1i}(t_2) = \frac{r_{1i}}{2}$ and $t_1 <t_2$ be the last time before $t_2$ such that $d_{1i}(t_1) = r_{1i}$. On $[t_1,t_2]$, $0<d_{1i}\leq r_{1i}$. Since $d_{1i}$ is absolutely continuous
    \begin{align}
        d_{1i}(t_2) - d_{1i}(t_1) = -\frac{r_{1i}}{2} = \int_{t_1}^{t_2} \dot{d}_{1i}(\tau) d\tau \geq -M_{1i}(t_2 - t_1)
    \end{align}
    which implies
    \begin{equation}
        t_2 - t_1 \geq \delta t_{1i} = \frac{r_{1i}}{2M_{1i}} > 0
    \end{equation}
    Inside $(\frac{r_{1i}}{2},r_{1i})$, using $[a]^+_\lambda \geq a$, the update law for $\dot{\lambda}_{1i}$ is bounded by
    \begin{equation}
        \dot{\lambda}_{1i}  \geq   -\alpha\lambda_{1i} + (\Gamma_1^{-1})_{ii}
\underline{c}_{1i}  
\end{equation}
where the last inequality is developed using $c(\hat{\theta}) \geq \underline{c}_{1i} >0$ when $0 < d_{1i}\leq r_{1i}$. Using comparison lemma, we have
\begin{align}
    \lambda_{1i}(t) \geq \lambda_{1i}(t_1)e^{-\alpha(t-t_1)} + (\Gamma_1^{-1})_{ii}\frac{\underline{c}_{1i}}{\alpha}(1 - e^{-\alpha(t-t_1)})
\end{align}
Since $\lambda_{1i}(t_1) \geq 0$,
\begin{equation}
    \lambda_{1i}(t_2) \geq (\Gamma_1^{-1})_{ii}\frac{\underline{c}_{1i}}{\alpha}(1 - e^{-\alpha\delta t_{1i}}) = \underline{\lambda}_{li} >0
\end{equation}
As long as $d_{1i} \leq r_{1i}$,, for $t>t_2$, $\lambda_{1i}(t) \geq \underline{\lambda}_{1i}$.

Since $\beta_{1i}(d_{1i})\rightarrow \infty$ as $d_{1i} \rightarrow 0^+$, $\exists\epsilon_{1i}$, $0 < \epsilon_{1i} \leq \frac{r_{1i}}{2}$ such that
\begin{equation}
    p_{ii}\underline{\lambda}_{1i}\beta_{1i}(\epsilon_{1i}) > M_{1i}
\end{equation}
which using (\ref{eq:dDotLowerBound}) and $\beta_{1i}(d_{1i}) \geq \beta_{1i}(\epsilon_{1i})$ for $0<d_{1i} \leq \epsilon_{1i}$ from (\ref{eq:betaDef}), implies for every $0<d_{1i} \leq \epsilon_{1i}$,
\begin{equation}
    \dot{d}_{1i} \geq \bar{d}_{1i} = p_{ii}\underline{\lambda}_{1i}\beta_{1i}(\epsilon_{1i}) - M_{1i} > 0
\end{equation}
Assume that, for contradiction, there exists $t^+ < \tau_\mathrm{max}$ such that $d_{1i}(t^+) = \epsilon_{1i}$. and suppose there exists $t > t^+$, such that $0< d_{1i}(\tau) \leq \epsilon_{1i}$, for $\tau \in [t^+,t]$ with $d_{1i}(t) < \epsilon_{1i}$. On any time interval $[t^+,t]$, due to absolute continuity, over which $d_{1i} \leq \epsilon_{1i}$, we have
\begin{equation}
    d_{1i}(t) - d_{1i}(t^+) = \int_{t^+}^t \dot{d}_{1i}(\tau) d\tau \geq \bar{d}_{1i}(t-t^+) > 0
\end{equation}
which implies that $d_{1i}(t)>\epsilon_{1i}$, contradicting the assumption that $d_{1i}<\epsilon_{1i}$. Hence,
\begin{equation}
    d_{1i}(t) \geq \epsilon_{1i} > 0, \quad \forall t < \tau_{\max}
\end{equation}

\paragraph{Behavior of $d_{2i}$ close to $\bar{\theta}_i$} The proof to show that $d_{2i} \geq \epsilon_{2i} > 0, \; \forall t < \tau_{\mathrm
max}$ follows along similar lines. Consider an interval $0<r_{2i} < \mathrm{min}\{d_{2i}(t_0),\delta_{2i},\frac{w_i}{2} \}$. On the interval $0 < d_{2i} < r_{2i}$, the terms $q_i + p_{ii}\lambda_{1i}\beta_{1i}(d_{1i})$ are bounded and we obtain
\begin{equation}
    \dot{d}_{2i} \geq -M_{2i} + p_{ii}\lambda_{2i} \beta_{2i}(d_{2i})
\end{equation}
where $M_{2i}$ is a positive constant. Using the analogous contradiction and comparison argument for $\lambda_{2i}$ gives a constant $\epsilon_{2i}>0$ such that
\begin{equation}
    d_{2i}(t) \geq \epsilon_{2i} > 0, \quad \forall t < \tau_\mathrm{max}
\end{equation}
Using bounds on $d_{1i}$ and $d_{2i}$, we obtain
\begin{equation}
    \underline{\theta}_i + \epsilon_{1i} \leq \hat{\theta}_i(t) \leq \bar{\theta}_i - \epsilon_{2i}, \quad \forall t < \tau_{\mathrm{max}} \label{eq:CompatBoundThetaHat}
\end{equation}
which implies $\hat{\theta}$ must remain within the compact set $\Theta^\kappa \subset \Theta^c,\;\forall t <\tau_\mathrm{max}$.

Using (\ref{eq:Lemma1Bounds}) and (\ref{eq:CompatBoundThetaHat}), the maximal solution satisfies $\zeta(t) \in K$ for $t \in [t_0,\tau_\mathrm{max})$, where $K:=\{e:\Vert e \Vert \leq e_m \} \times \Theta^\kappa \times \{ \lambda \in \mathbb{R}_{\geq 0}^{2p}: \Vert \lambda \Vert \leq \lambda_m \}$ is compact. Suppose, for contradiction, that $\tau_{\max}<\infty$. Since $\hat{\theta} \in \Theta^\kappa$, $c_j(\hat{\theta}(t))$ and $\nabla c_j(\hat{\theta}(t))$ are uniformly bounded and locally Lipschitz on the neighbor of  $\Theta^\kappa$. Together with (\ref{eq:Lemma1Bounds}) and the boundedness of $Y(x)$ and $S(t)$, the unprojected vector field (\ref{eq:ErrorDyn}) and (\ref{eq:ParameterUpdateLaw}) is locally bounded on $[t_0,\tau_\mathrm{\max}) \times K$. Local existence results of time varying projected dynamical systems (\cite{heemels2023existence}, Corollary 2) yield a Caratheodory solution starting from every point of $K$, hence, a maximal solution can be continued beyond $\tau_\mathrm{max}$, contradicting maximality. Therefore, $\tau_{\max}=\infty$. Consequently, $\Theta^c$ is forward invariant.
\end{proof}
Stability of the error system (\ref{eq:ErrorDyn}), (\ref{eq:thetaTildeDot}) and (\ref{eq:ParameterUpdateLaw}) is now analyzed using Lyapunov analysis. 
\begin{theorem}
    Provided Assumptions \ref{ass:MarginParam}-\ref{ass:HistoryStackRegularity} are satisfied, $\zeta(t_0) \in \mathcal{D}$. Then, for the system in (\ref{eq:SystemModel}), the constrained parameter update law (\ref{eq:ParameterUpdateLaw}) and the adaptive controller (\ref{eq:Control}) ensure that $\hat{\theta} \in \Theta^c$ and $e(t)$ and $\tilde{\theta}$  are uniformly ultimately bounded (UUB).
\end{theorem}
\begin{proof}
From Lemma 2, $\tau_{\max}=\infty$ and $\Theta_c$ is forward invariant. 
From Lemma 1, for $t \geq t_H$,
\begin{equation}
V(t) \leq V(t_H)e^{-\beta(t-t_H)} +\frac{D_x}{\beta} \left(1-e^{-\beta(t-t_H)}\right).
\end{equation}
Therefore, $\limsup_{t\to\infty} V(t) \leq \frac{D_x}{\beta}$. Using the lower bound in (24),
$\Lambda_{\min}\|z(t)\|^2 \leq V(t)$, gives
\begin{equation}
\limsup_{t\to\infty}\|z(t)\| \leq \sqrt{\frac{D_x}{\beta\Lambda_{\min}}}.
\end{equation}
Since $\limsup_{t\to\infty}\Vert[e^T,\tilde{\theta}^T]\Vert \leq \limsup_{t\to\infty}\Vert z(t) \Vert \leq \sqrt{\frac{D_x}{\beta\Lambda_{\min}}}$, $e(t)$ and $\tilde{\theta}$ are UUB.
\end{proof}

\section{Simulations}
Simulations are carried out to demonstrate the applicability of the constrained parameter update law in adaptive control setting. The following system dynamics similar to \cite{parikh2019integral} is considered
\begin{align}
    \dot{x} = \left[ \begin{matrix}
        x_1^2 & sin(x_2) & 0 & 0 \\
        0 & x_2sin(x_1) & x_1 & x_1x_2
    \end{matrix}\right]\theta + u
\end{align}
where $x = [x_1\;x_2]^T \in \mathbb{R}^2$ is the state, $u\in\mathbb{R}^2$ is the control input, $\theta \in \mathbb{R}^4$ are the true parameters whose values are given by $\theta = [5 \; 10 \; 15 \; 20]^T$.
The desired trajectory is given by
\begin{equation}
    x_d(t) = 10(1-e^{-0.1t})\left[ \begin{matrix}
    \mathrm{sin}(2t) \\ 0.4\mathrm{cos}(3t) 
    \end{matrix}\right]
\end{equation}

\subsection{Simulation with Bounds on $\hat{\theta}_i$ using Inverse Barrier}

The controller in (\ref{eq:Control}) is implemented along with the constrained parameter update law in (\ref{eq:ParameterUpdateLaw}). The system state is initialized to $x(t_0) = [10\;5]^T$, the parameters are initialized to $\hat{\theta}(t_0) = [4.5\;8\;12\;15]^T$ which are within the parameter bounds, the dynamic Barrier weights for upper and lower bound constraints are initialized to $\lambda_1(t_0)=\lambda_2(t_0) = [5\;5\;5\;5]^T$. The control gain is selected as $k=10\mathbb{I}_{2\times 2}$ and for the parameter adaptation law $k_{cl} = 0.1$, $P = \mathrm{diag}(0.025,0.04,0.12,0.008)$, $\Gamma^{-1} = 0.15\mathbb{I}_{4 \times 4}$ and $\alpha = 0.1$. For the shifted-scaled inverse Barriers margin parameters are selected as $\delta_1 = \delta_2 = 0.5\mathbf{1}_4$. The upper and lower bounds on parameters are selected as $\bar{\theta} = [6,\;12,\;17,\;22]^T$ and $\underline{\theta} = [3,\;6,\;10,\;12]^T$.

\begin{figure}
    \centering
    \includegraphics[width=\linewidth]{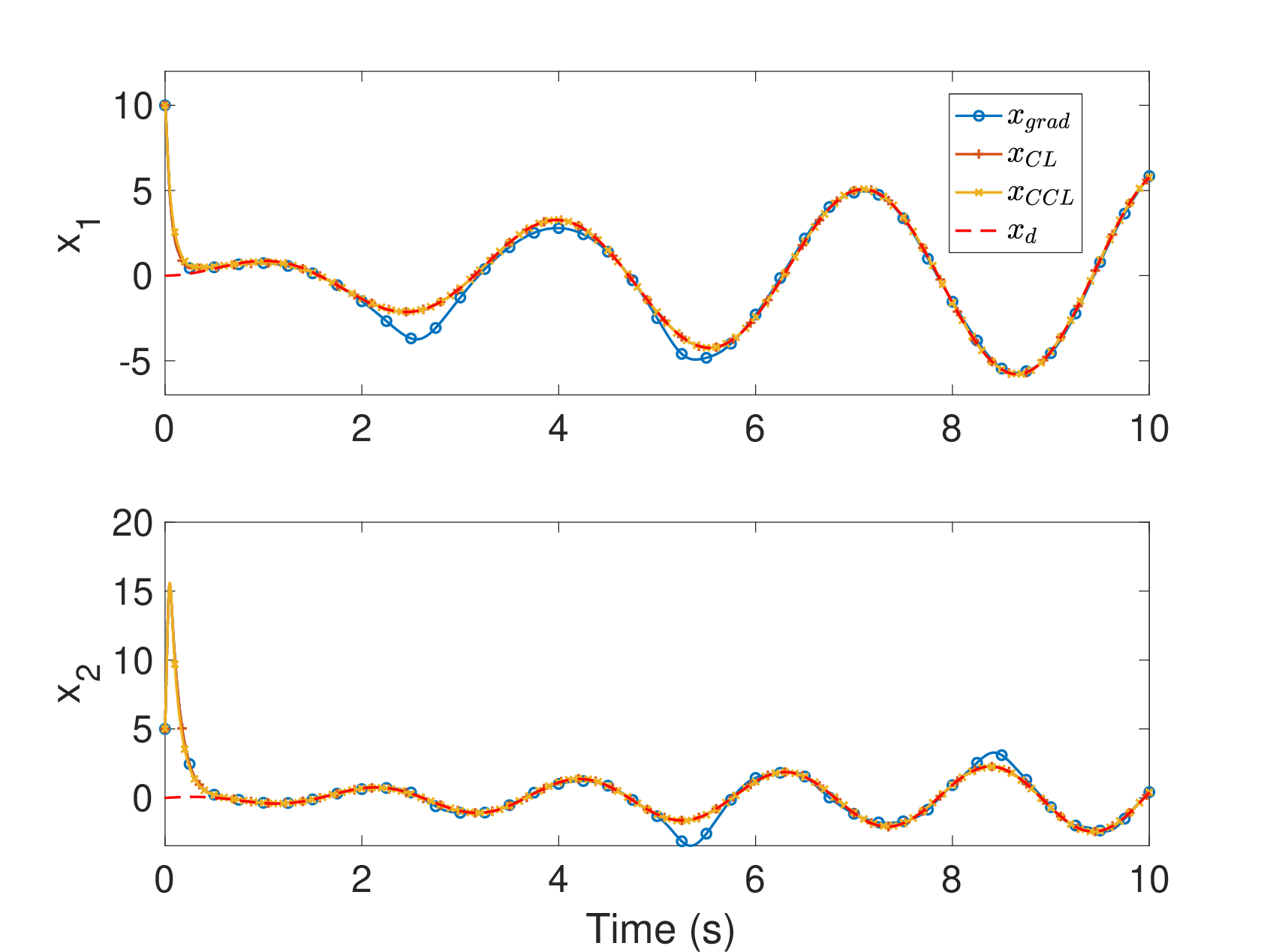}
    \caption{Trajectory tracking with different parameter update laws.}
    \label{fig:TrajTracking}
\end{figure}

\begin{figure}
    \centering
    \includegraphics[width=1.0\linewidth]{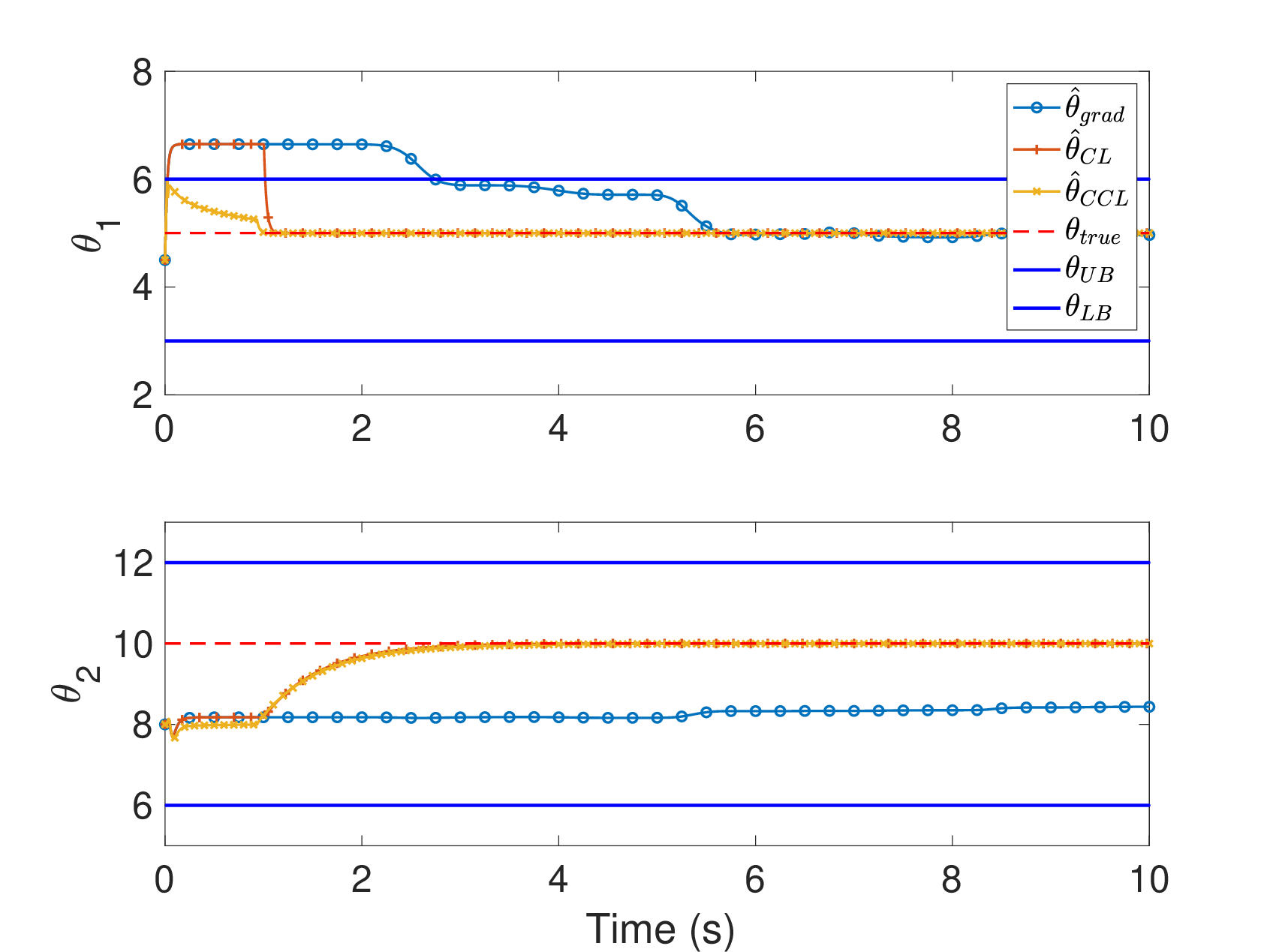}
    \caption{Parameter convergence with different parameter update laws for $\theta_1$ and $\theta_2$.}
    \label{fig:ParamComparison1}
\end{figure}

\begin{figure}
    \centering
    \includegraphics[width=1.0\linewidth]{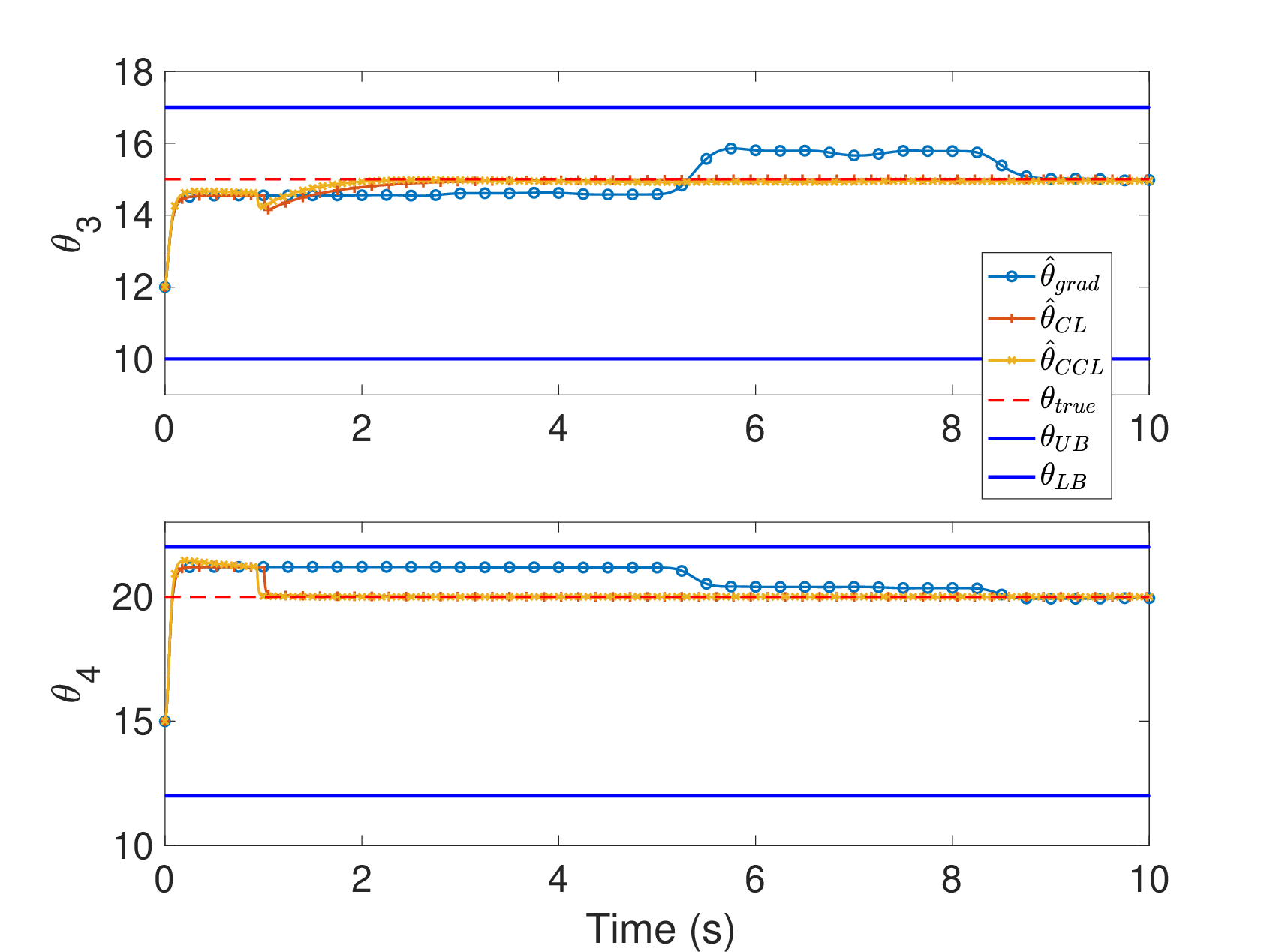}
    \caption{Parameter convergence with different parameter update laws for $\theta_3$ and $\theta_4$.}
    \label{fig:ParamComparison2}
\end{figure}

\begin{figure}
    \centering
    \includegraphics[width=\linewidth]{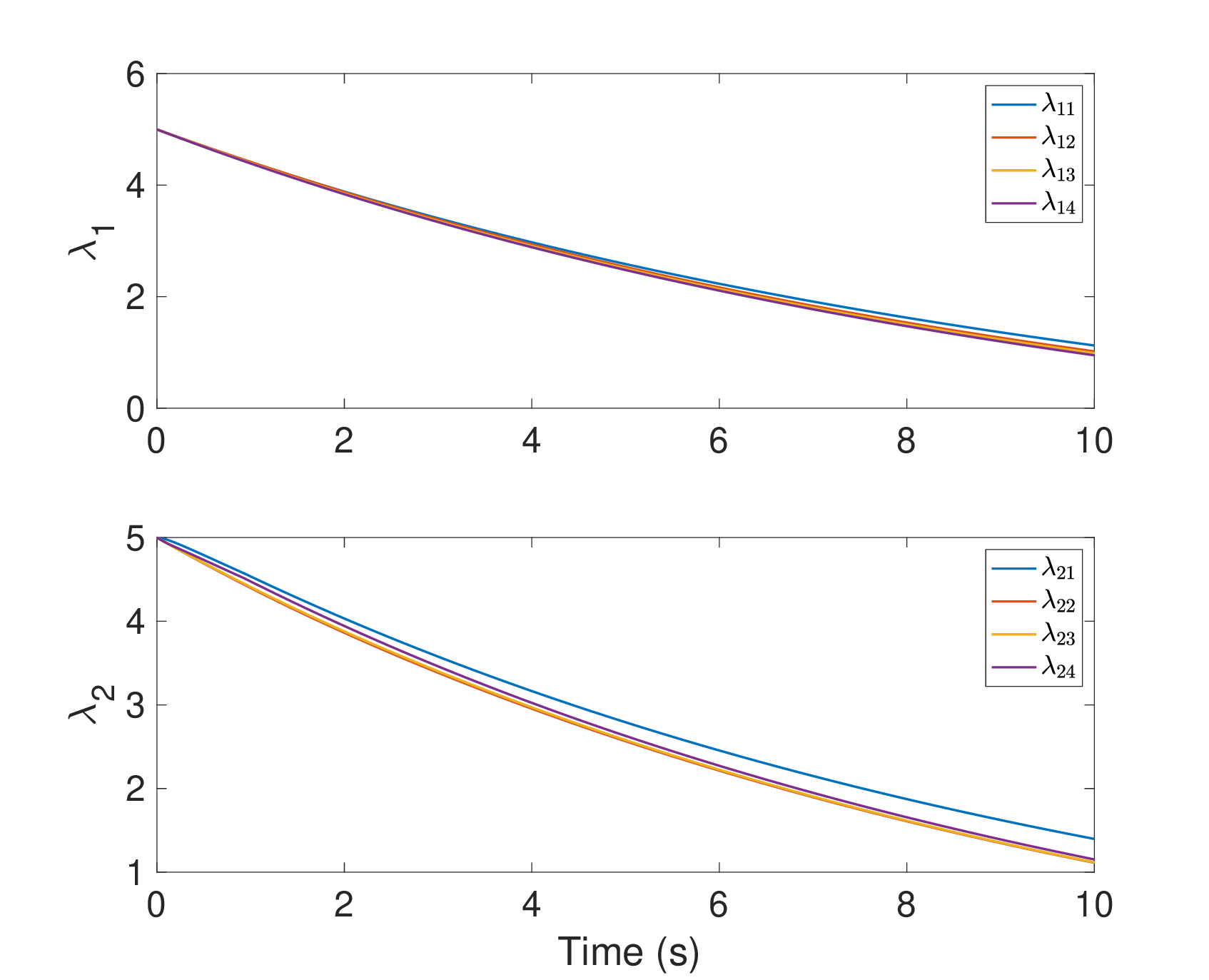}
    \caption{Evolution of dynamic Barrier weights for inverse Barrier constraint.}
    \label{fig:Lagrange1}
\end{figure}

The performance of the adaptive controller with constrained parameter update law is compared to that with gradient-based parameter update law and CL-based parameter update law for its ability to keep the parameter estimates bounded. The results are shown in Figs. \ref{fig:TrajTracking}, \ref{fig:ParamComparison1}-\ref{fig:ParamComparison2} and \ref{fig:Lagrange1}. It is observed from Figs. \ref{fig:ParamComparison1}-\ref{fig:ParamComparison2} that the parameters estimated by the proposed constrained parameter update law remain strictly within the prescribed bounds and the parameter convergence improves compared to purely gradient-based update law. Compared to the CL-based parameter estimation due to the presence of constraints in the parameter estimation law, the estimates stay strictly within the bounds and move toward true parameters during the transient as can be seen from Figs. \ref{fig:ParamComparison1}-\ref{fig:ParamComparison2}. The trajectory tracking performance is better compared with the gradient-based parameter update law and comparable when the parameters are estimated using CL-based parameter update law.



\subsection{Simulation with Bounds on $\hat{\theta}_i$ using log-Barrier}
\begin{figure}
    \centering
    \includegraphics[width=\linewidth]{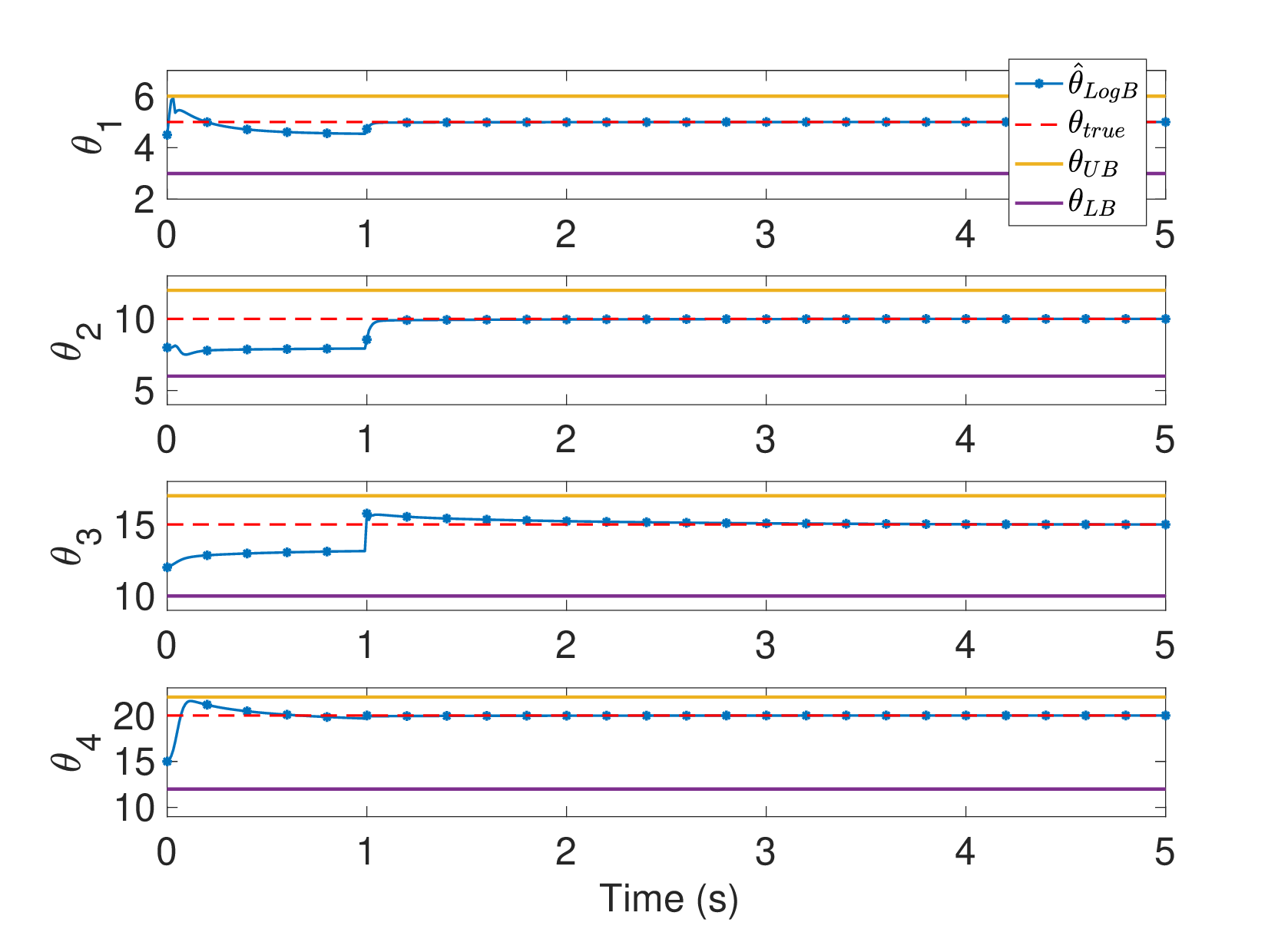}
    \caption{Parameter estimation using constrained parameter update law with log-Barrier constraints.}
    \label{fig:ParamLogB}
\end{figure}

\begin{figure}
    \centering
    \includegraphics[width=\linewidth]{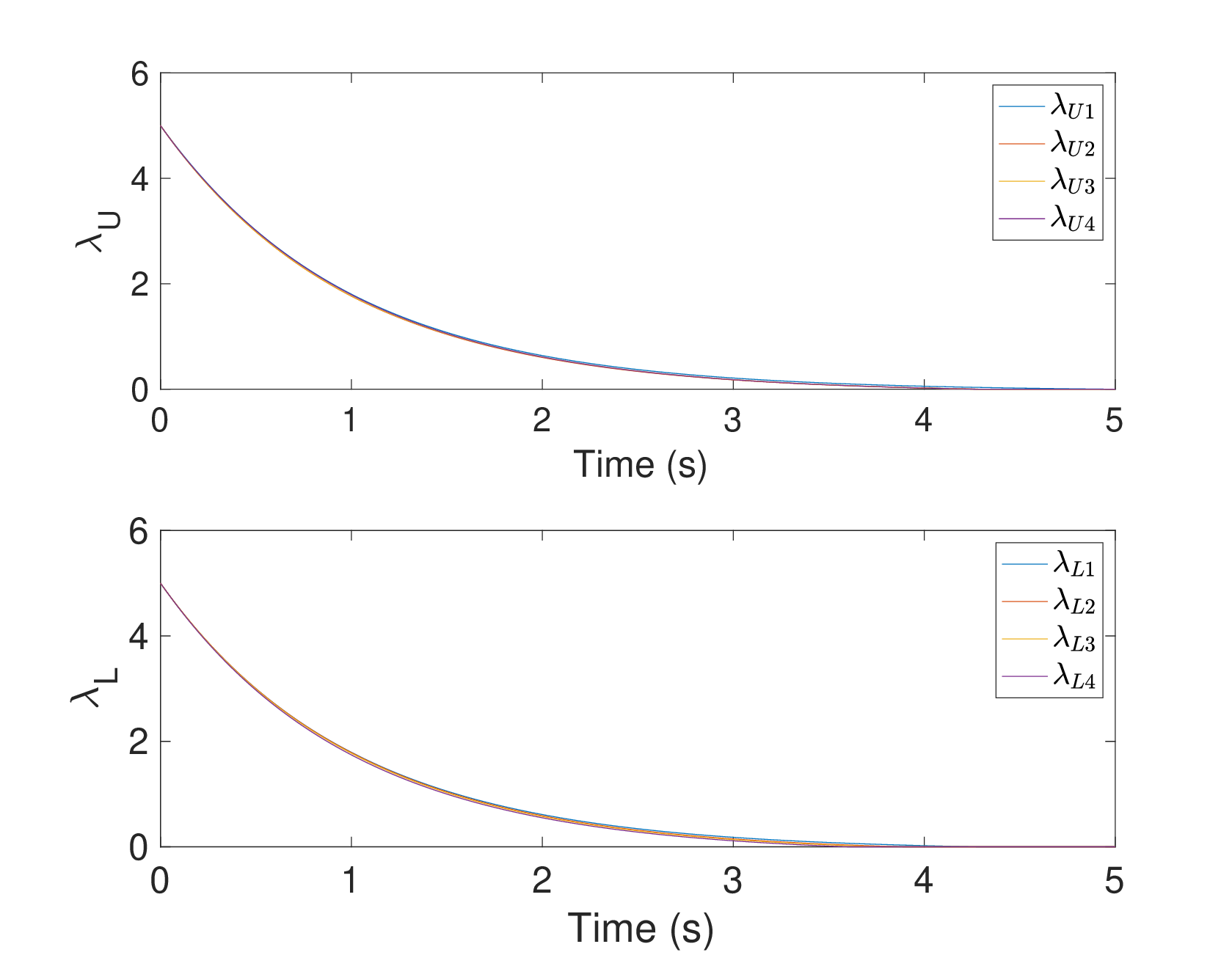}
    \caption{Evolution of the dynamic Barrier weights for log Barrier constraint.}
    \label{fig:LagrangeLogB}
\end{figure}

The controller is implemented using log-Barrier constraints on $\hat{\theta}$ for the same dynamics and desired trajectories as that of previous simulations. The scaled log-barrier constraints are formulated as using $\delta_{1i} = \delta_{2i} = 0.5, \forall i =\{1,..,4\}$. The parameter estimates and the dynamic Barrier weights are initialized to the same values as that of previous simulations. The control gain is selected as $k=10\mathbb{I}_{2 \times 2}$, and the parameter adaptation gains are selected as $k_{cl} = \mathrm{diag}(0.02,0.5,1.5,0.02)$, $P= 1.5 \mathrm{diag}(0.8,0.1,1,0.9)$, $\Gamma^{-1} = 0.05\mathbb{I}_{4 \times 4}$ and $\alpha = 1$. From the results shown in Figs. \ref{fig:ParamLogB}-\ref{fig:LagrangeLogB}, it is observed that the parameter estimates converge to their true values and parameter estimate bounds are maintained.

\section{Conclusion}
In this paper, a constrained parameter update law is developed for adaptive tracking control using a Barrier constraint formulation, where dynamic barrier weights are updated using an update law. The approach of deriving adaptive parameter update law using a minimization problem of an objective function allows for incorporation of constraints on the parameter update law design. Lyapunov stability analysis of the tracking error, parameter estimation error and barrier weight dynamics is conducted which shows that the tracking error and parameter estimates converge to a small bound around zero while the parameter estimates stay within a confined set. Simulation results demonstrate the benefits of the proposed approach in improving transient performance and enforcing parameter constraints.

\bibliographystyle{IEEEtran}
\bibliography{RCL_Complete}

\end{document}